\documentclass[12pt]{article}
\usepackage{amsmath,amsfonts,graphics,color,amssymb,amstext}
\newtheorem{thm}{Theorem} \newtheorem{lemma}[thm]{Lemma}
 
\newenvironment{pf} {\noindent{\sc Proof. }}{{\hfill
$\Box$}\par\vskip2\parsep}
\newenvironment{pfof}[1]
{\par\vskip2\parsep\noindent{\sc Proof of\ #1. }}{{\hfill $\Box$}
  \par\vskip2\parsep}

\newcommand{\py}{\mbox{\bf P}}
\newcommand{\ey}{\mbox{\bf E}}
\newcommand{\rd}{{\mathbb R}^d}

\newcommand{\dof}{\bf\boldmath}
\newcommand{\ind}{\mbox{\boldmath $1$}}

\newcommand{\leb}{{\mathcal L}}
\newcommand{\ld}{{\mathcal D}}
\newcommand{\hx}{\ensuremath{{\mathbb H^{*}}}}
\newcommand{\as}{a.s.\ }
\newcommand{\pch}{p_{\text{\rm c}}^{\text{\rm site}}({\mathbb H})}
\newcommand{\pc}{p_{\rm c}}
\newcommand{\bigmid}{\;\big|\;}
\newcounter{mycount}
\newenvironment{mylist}{\begin{list}{(\roman{mycount})}%
{\usecounter{mycount}\itemsep 0pt}}{\end{list}}
\usepackage{graphics}

\begin{document}

\title{A Percolating Hard Sphere Model}
\author{
Codina Cotar\thanks{
Funded in part by NSERC},
Alexander E. Holroyd\thanks{
Funded in part by an NSERC Discovery Grant, and by CPAM and MSRI.}
 \ and David Revelle\thanks{
Funded in part by an NSF Postdoctoral Fellowship.} }
\date{July 25, 2006}

\maketitle
\renewcommand{\thefootnote}{}
\renewcommand{\thefootnote}{\arabic{footnote}}

\begin{abstract}
Given a homogeneous Poisson point process in ${\mathbb R}^d$,
H\"aggstr{\"o}m and Meester \cite{hagmee} asked whether it is
possible to place spheres (of differing radii) centred at the
points, in a translation-invariant way, so that the spheres do not
overlap but there is an unbounded component of touching spheres.
We prove that the answer is yes in sufficiently high dimension.
\end{abstract}

\section{Introduction}

A {\dof sphere process} is a simple point process $\Lambda$ on
$\rd\times [0,\infty)$. The support of $\Lambda$ is the random set
$[\Lambda]=\{(x,r)\in\rd\times[0,\infty):\Lambda(\{(x,r)\})=1\}$. If
$(x,r)\in[\Lambda]$, we say that there is a sphere of radius
$r\in[0,\infty)$ at $x\in \rd$. The {\dof centre process}
$\widetilde{\Lambda}$ is the point process on $\rd$ given by the
projection
$\widetilde{\Lambda}(\cdot)=\Lambda(\cdot\times[0,\infty))$.  We say
that $\Lambda$ is a {\dof Poisson sphere process} if
$\widetilde{\Lambda}$ is a homogeneous Poisson process.

A {\dof hard sphere process} is one in which the interiors of the
spheres do not overlap;
that is, almost surely
$$\rho(x,y)\geq r+s \text{ for any distinct }
(x,r),(y,s)\in[\Lambda],$$
where $\rho$ denotes Euclidean distance.
(Note in particular that, since all the radii are non-negative,
in a hard sphere process $\Lambda$ no sphere may contain any points of
$\widetilde{\Lambda}$ other than its own centre.)
For $z\in\rd$ the shifted sphere process $\Lambda+z$ is defined by
$(\Lambda+z)(A)=\Lambda(A-z)$, where $A-z:=\{(x-z,r):(x,r)\in
A\}$.    The sphere process
$\Lambda$ is (translation-){\dof invariant} if $\Lambda+z$ and $\Lambda$
are equal in law for all $z\in\rd$.

Let $G(\Lambda)\subseteq\rd$ be the random set covered by all the
spheres:
$$G(\Lambda)=\big\{y\in\rd: \rho(y,x)\leq r \text{ for some }
(x,r)\in[\Lambda]\big\}.$$ The connected components of
$G(\Lambda)$ are called {\dof clusters}. We say that $\Lambda$
{\dof percolates} if there is an unbounded cluster. Our main
result is the following.
\begin{thm}
\label{main1}
For all $d\ge 45$, there exists an invariant Poisson hard
sphere process $\Lambda$ which percolates almost surely.
\end{thm}

H\"aggstr{\"o}m and Meester \cite{hagmee} studied several invariant
continuum percolation processes, and proved that the ``dynamic
lily-pond model'' does not percolate.  This is a Poisson hard sphere
process in which spheres grow from all Poisson points at the same rate,
and whenever two spheres touch, they both stop growing.  H\"aggstr{\"o}m
and Meester asked whether there exists an invariant Poisson hard
sphere process which percolates.

In dimension $d=1$ it is easy to see that any Poisson hard sphere
process almost surely does not percolate.  For consider $4$ consecutive
points $x_1<x_2<x_3<x_4$ of the Poisson process which satisfy
$x_3-x_2>(x_2-x_1)+(x_4-x_3)$; then the spheres centred at $x_2$ and
$x_3$ cannot touch.  Since such a configuration of $4$ points appears
infinitely often in a Poisson process, no percolation is possible.  It
is unknown whether there exists a percolating Poisson hard sphere
process in dimensions $2\leq d\leq 44$ (even without the requirement
of invariance).

Our explicit bound of $45$ in Theorem \ref{main1} could probably
be reduced with some effort.  Our construction certainly cannot be
adapted to work in dimension less than $3$
(and probably not in dimension that low).

Jonasson \cite{jonasson} showed that in the hyperbolic plane $H^2$
there exists a hard sphere process that percolates when the
Poisson process has sufficiently low intensity, but was unable to
determine what happens for high intensity.  (In $\rd$, scaling shows that
the intensity of the process is immaterial).

Note that our definition of a hard sphere process allows spheres of
radius zero.  Our proof of Theorem~\ref{main1} can in fact be adapted
to prove the existence of an invariant percolating Poisson hard sphere
process in which the spheres all have positive radii.  We explain this
in a remark at the end of Section \ref{average}.

Our proof of Theorem~\ref{main1} is in two parts.  First we construct
a non-invariant hard sphere process $\Gamma$ which percolates.  Then
we convert this to an invariant process $\Lambda$ by
``stationarizing'' - applying a uniform random translation in a large
ball, and taking a limit.  The non-invariant construction of $\Gamma$
proceeds as follows.  Starting from a Poisson process, we attempt to
grow an unbounded cluster iteratively.  At each step we try to choose
a radius for a new sphere (centred at some Poisson point) so that it
touches the cluster constructed so far.  This will be possible only if
the proposed new sphere contains no other Poisson points, and we need
to ensure this happens sufficiently often that the cluster can
continue to grow.  We do this by comparison with a certain
two-dimensional percolation process; in sufficiently high dimension we
can arrange that the probability of success at each step exceeds the
relevant critical probability.

To ensure that the stationarized version $\Lambda$ also
percolates, the un\-bound\-ed clusters of $\Gamma$ should occupy a
positive fraction of space -- see below for a precise statement.
We will achieve this by repeating the essentially two-dimensional
construction described above in infinitely many ``layers''
throughout $\rd$.

Let $\leb$ denote Lebesgue measure on $\rd$, and denote the
Euclidean ball $B(x,r)=B_d(x,r):=\{y\in\rd: \rho(x,y)<r\}$, and
the origin $0=(0,\ldots,0)$.
 Define the {\dof lower density} of a set
$A\subseteq\rd$ to be
$$\ld(A) = \liminf_{r\to\infty} \frac{\leb\big(A\cap B(0,r)\big)}{\leb B(0,r)}.$$
For $s>0$ define the {\dof $s$-neighbourhood} of a set $A\subseteq\rd$ to be
$$A^{\{s\}}=\bigcup_{x\in A} B(x,s).$$
Using
the construction sketched above, we shall prove the following
result, from which Theorem \ref{main1} will be deduced.
\begin{thm}
\label{noninv} For $d\geq 45$ there exists a (not necessarily
invariant) Poisson hard sphere process $\Gamma$ which percolates
almost surely. Furthermore, $\Gamma$ can be chosen so as to have
the following additional properties.
\begin{mylist}
\item The union $I$ of all unbounded clusters satisfies
\hbox{$\lim_{a\to\infty} \ey\ld(I^{\{a\}})=1$.}
\item
There exists a constant $K(d)<\infty$ such that $\Gamma(\rd\times
(K,\infty))=0$ almost surely; that is, there are no spheres of radius
greater than $K$.
\end{mylist}
\end{thm}

The article is organized as follows.  In Section
\ref{construction} we give the construction that gives rise to the
non-invariant process in Theorem \ref{noninv}.  In Section
\ref{proof} we prove that the resulting process does indeed have
the properties stated in Theorem \ref{noninv}.  In Section
\ref{average} we deduce Theorem \ref{main1}.

\section{Construction}
\label{construction}

In this section we describe the construction of the hard sphere
process in Theorem \ref{noninv}.  We will prove that it has the
required properties in the next section.  The construction is
given in terms of parameters $\mu=0.75$, $\delta=0.1$,
$\epsilon=0.01$ and $C$, where $C=C(d)$ is a (large) constant to
be chosen later; the values for these parameters are chosen so
that the construction yields an unbounded cluster.  Let
$\lambda=\lambda(d)>0$ be another constant to be chosen later, and
let $\Pi$ be a homogeneous Poisson process of intensity $\lambda$,
and denote its support $[\Pi]$.
We will construct a hard sphere process $\Gamma$ such that
$\widetilde{\Gamma}=\Pi$.

\begin{figure}
\begin{center}
\resizebox{4cm}{!}{\includegraphics{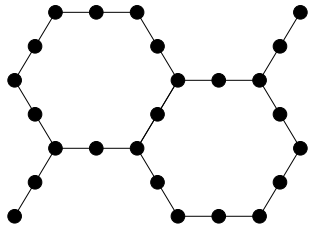}}
\end{center}
\caption{Part of the lattice \hx \label{hex}}
\end{figure}
Let ${\mathbb H}$ be the hexagonal lattice in ${\mathbb R}^2$,
with edge-length $2$.  (Thus, the faces are regular hexagons of
side $2$, and there is a vertex at the origin $0=(0,0)$, say.) Let
${\mathbb H^{*}}$ be the graph formed from ${\mathbb H}$ by adding
an extra vertex in the middle of each edge; see Figure \ref{hex}.
We will call the vertices of the original lattice ${\mathbb H}$
\textbf{site vertices} and the extra vertices of ${\mathbb H^{*}}$
\textbf{bond vertices}.  Also fix some arbitrary well-ordering of
the vertex set $V(\hx)$ of \hx.

\begin{figure}
\begin{center}
\resizebox{6cm}{!}{\includegraphics{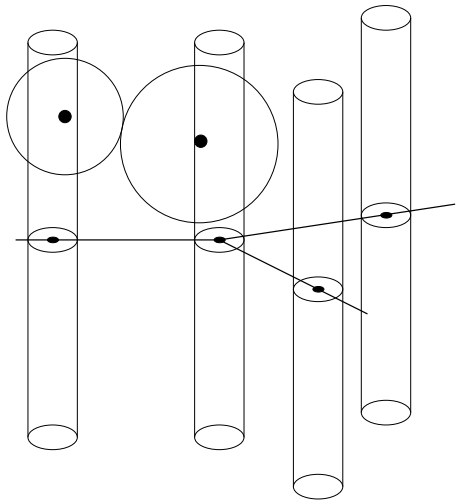}}
\end{center}
\caption{An example of the construction in dimension $d=3$.
Shown are the cells centred at the origin $0$ and at the three
neighbouring bond-vertices $u,v,w$, together with two touching balls
$B(x_0,r_0)$ and $B(x_u,r_u)$.
 \label{cells}}
\end{figure}

Fix the dimension $d\geq 3$.  For any vertex $v\in V(\hx)$, define
the {\dof cell} centred at $v$ to be the set
$$W(v):=B_2(v,\epsilon)\times B_{d-2}(0,C) \subset\rd.$$
See Figure \ref{cells}.

First, here is a brief description of the construction.  We will
attempt to place spheres with their centres in distinct cells in
such a way that spheres at adjacent cells touch; see Figure
\ref{cells}. All the spheres will have radii in
$[\mu-\delta,\mu+\delta]$. Note that any two non-adjacent vertices
of \hx\ are at distance at least $\surd 3$. Since
$2(\mu+\delta+\epsilon)=1.72<\surd 3$, it will be impossible for
spheres centred in non-adjacent cells to touch or overlap.

We will explore the lattice iteratively starting from the origin,
attempting to construct an unbounded cluster.  Each step of this
exploration will have two parts.  Firstly, for some vertex $w\in
V(\hx)$ we try to find a point $y_w\in[\Pi]\cap W(w)$ such that it
is possible to place a sphere of some radius $r_w\in
[\mu-\delta,\mu+\delta]$ at $y_w$ that touches one of the existing
spheres and does not overlap any of the existing spheres.
Secondly, we check to see whether or not there are any other
points of the Poisson process within $B(y_w, r_w)$. If there is
none, our construction succeeds and we let $x_w=y_w$. If either
part fails, we let $x_w=\Delta.$

At each step of the construction, each vertex $v\in V(\hx)$ will be
either \textbf{good} (meaning we succeeded in constructing a sphere
$B(x_v,r_v)$ at $v$), \textbf{bad} (meaning the construction failed),
or {\dof unexplored}.

Using \hx\ (instead of ${\mathbb H}$, say) will enable us to
explore in such a way that every time we explore a new vertex it
is adjacent to exactly one good vertex and no bad vertices.  This
will simplify our arguments.  Other choices of lattice are
possible, such as the square lattice in place of ${\mathbb H}$ --
see the remark at the end of Section \ref{proof}.

Here is a formal description of the construction.  As remarked
earlier, we will try to construct a cluster in each of a set of
two-dimensional ``layers''.

\subsection*{First Layer}

Start with all vertices in $V(\hx)$ unexplored.  We perform a
sequence of steps $0,1,2,\ldots$.  One new vertex (or sometimes
two) will be explored at each step.

\subsubsection*{Step $\mathbf{0}$}

Let $w_0=0$; we will start our exploration at the site vertex
$w_0$.  Temporarily denote its cell $W=W(w_0)$. Consider two
cases:
\begin{itemize}
\item[] \begin{itemize} \item[\emph{Case 1}:] If $W\cap
[\Pi]=\emptyset$, then let $x_0=y_0=\Delta.$ \item[\emph{Case 2}:]
If $W\cap [\Pi]\not=\emptyset$,  then pick $y_0$ uniformly at
random from the set $W\cap [\Pi]$ (conditional on $W\cap [\Pi]$).
Take $r_0=\mu$.  Now take
$$x_0=
\left\{
\begin{array}{lcc}
y_0& \mbox{if} & B(y_0, r_0)\cap [\Pi]=\{y_0\}; \\
\Delta & \mbox{if} & B(y_0, r_0)\cap [\Pi]\supsetneq\{y_0\}.
\end{array}
\right.$$
\end{itemize}
\end{itemize}
Declare $w_0=0$ to be bad if $x_0=\Delta$, and good otherwise.

\subsubsection*{Step $\mathbf{n}$ ($\mathbf{n\geq 1}$)}

Step $n$ consists of parts (a)--(c).

\begin{enumerate}

\item [(a)]  We first choose an unexplored vertex $w_n$ to
explore, according to the following rules.
\begin{mylist}
\item If there is an unexplored site vertex $w$ adjacent to some
good bond vertex and two unexplored bond vertices, choose $w_n$ to
be the first such $w$ in the ordering on $V(\hx)$. \item If there
is no $w$ as in (i), but there is an unexplored bond vertex $w$
adjacent to a good site vertex and an unexplored site vertex,
choose $w_n$ to be the first such $w$. \item If neither (i) nor
(ii) hold then \textbf{stop}, and proceed to `Subsequent Layers'
below.
\end{mylist}

\item[(b)] Temporarily write $w=w_n$ for the vertex chosen in (a).
The rules in (a) ensure that $w$ is unexplored, and it has exactly
one good neighbour, $v$ say, while all the other neighbours of $w$
are unexplored.  We will try to construct $x_w\in W(w)$,
the centre of a sphere tangent to $B(x_v,r_v)$.

Define
$$U(v)=B(x_v, r_v+\mu+\delta)\setminus B(x_v,r_v+\mu-\delta);$$
this is the set of possible centres for a sphere of radius $\in
[\mu-\delta,\mu+\delta]$ that touches $B(x_v,r_v)$.  Let
$$S=W(w)\cap U(v).$$
Now consider two cases:

\begin{itemize}
\item[] \begin{itemize}
\item[\emph{Case 1}:]
If $S\cap [\Pi]=\emptyset,$ then let $x_w=y_w=\Delta$.

\item[\emph{Case 2}:]
If $S\cap [\Pi]\not=\emptyset$, then pick $y_w$ uniformly at random from
$S\cap [\Pi]$ (conditional on $S\cap [\Pi]$), and let
$$r_w=\rho(x_v,y_w)-r_v$$
(by the definition of $V(v)$ we have $r_w\in[\mu-\delta,\mu+\delta]$).
Now take
$$x_w=
\left\{
\begin{array}{lcc}
y_w& \mbox{if} & B(y_w, r_w)\cap [\Pi]=\{y_w\}; \\
\Delta & \mbox{if} & B(y_w, r_w)\cap [\Pi]\supsetneq\{y_w\}.
\end{array}
\right.$$
\end{itemize}
\end{itemize}

\item [(c)] Declare the vertex $w$ to be bad if $x_w=\Delta$, and
good otherwise.

{\em In addition, if $w$ is a bond vertex and is
declared bad, then declare its remaining unexplored site
vertex neighbour to be bad also.}
\end{enumerate}

Continue either for an infinite sequence of steps,
or until we {\bf stop} in (a)(iii) above.

\subsection*{Subsequent Layers}

All the points $x_w$ constructed above (for good vertices $w$) lie
in the ``layer'' ${\mathbb R}^2\times B_{d-2}(0,C)$, and the radii
$r_w$ are at most $\mu+\delta$.  We want to repeat the
construction in other layers, ensuring that the spheres in
different layers do not overlap.  Therefore let
$L=2(C+\mu+\delta+1)$, and consider the lattice of points
$L\mathbb{Z}^{d-2}\subset{\mathbb R}^{d-2}$.  For each $z\in
L\mathbb{Z}^{d-2}$, repeat the ``first layer'' construction above,
but now in the layer ${\mathbb R}^2\times B_{d-2}(z,C)$ (using the
Poisson points in cells of the form $B_2(v,\epsilon)\times
B_{d-2}(z,C)$). This results in independent identically
distributed clusters in each of the layers.

\subsection*{Definition of $\Gamma$}

Finally, construct the hard sphere process $\Gamma$ by placing a
sphere centred at $x_w$ with radius $r_w$, for each good vertex
$w$, in every layer.  Also place a sphere of radius zero centred
at each remaining Poisson point $z\in[\Pi]$.

\section{Proof of Construction}
\label{proof}

In this section we prove Theorem \ref{noninv}, using the
construction of Section \ref{construction}.  We start by
assembling some tools.

\subsection*{Percolation}

Let ${\mathbb H}$ be the hexagonal lattice in ${\mathbb R}^2$ with
side-length 2.  Consider Bernoulli site percolation with parameter $p$
on ${\mathbb H}$. That is, each vertex is open with probability $p$
and closed otherwise, independently for different vertices. There
exists a critical probability $\pc=\pch<1$, with the property that if
$p>\pc$ then there is \as a unique infinite connected cluster of open
vertices (see \cite{grim}).  It is proved in \cite{wierman} that
$$\pc<0.794.$$

Now suppose $p>\pc$ and let $K\subset {\mathbb R}^2$ be the vertex
set of the infinite open cluster if contains $0$, and let
$K=\emptyset$ otherwise. Suppose $d\geq 3$ and let
$(K_z)_{z\in{\mathbb Z}^{d-2}}$ be a family of i.i.d.\ random sets
each with the same law as $K$. Fix any $L>0$, and define the
random set
$$Y=\bigcup_{z\in{\mathbb Z}^{d-2}} K_z \times \{Lz\}  \subset \rd.$$

\begin{lemma}
\label{perc-dens} For $p>\pch$ and any $L>0$, the set $Y$ defined
above satisfies $\ey\ld(Y^{\{a\}})\to 1$ as $a\to\infty$.
\end{lemma}

\begin{pf}
Let $I$ be the infinite open cluster, and let $\theta=\py(0\in
I)>0$. Define a random set $\widetilde{K}$ as follows. Flip a
$\theta$-coin independently of $I$, and let $\widetilde{K}=I$ with
probability $\theta$, and $\widetilde{K}=\emptyset$ with
probability $1-\theta$. By the Harris-FKG inequality (see for
example \cite{grim}) we see that $K$ stochastically dominates
$\widetilde{K}$.  Now let $(\widetilde{K}_z)_{z\in{\mathbb
Z}^{d-2}}$ be a family of i.i.d.\ random sets with the same law as
$\widetilde{K}$, and let $\widetilde{Y}=\bigcup_{z\in{\mathbb
Z}^{d-2}} \widetilde{K}_z \times \{Lz\} \subset \rd$. Then $Y$
dominates $\widetilde{Y}$, so $\ld(Y^{\{a\}})$ dominates
$\ld(\widetilde{Y}^{\{a\}})$. The random set
$\widetilde{Y}^{\{a\}}$ is invariant in law under isometries of
${\mathbb H}\times {\mathbb Z}^{d-2}$, so Fubini's theorem implies
$$
\ey\ld\big(\widetilde{Y}^{\{a\}}\big)
\geq \inf_{x\in\rd} \py\big( x\in\widetilde{Y}^{\{a\}} \big)
= \inf_{x\in\rd}\py\big(B(x,a)\cap\widetilde{Y}\neq\emptyset\big)
\to \py\big(\widetilde{Y}\neq\emptyset\big)
=1$$
as $a\to\infty$.
\end{pf}

\subsection*{Random points}

Let $\Pi$ be a point process in $\rd$. For $r>0$ we call a point
$x\in[\Pi]$ {\dof $r$-isolated} if there is no other point of
$[\Pi]$ within distance $r$ of $x$.

\begin{lemma}
\label{isolated} Let $\Pi$ be a homogeneous Poisson point process
of intensity $\lambda$ in $\rd$. Let $S\subseteq\rd$ be a Borel
set with $\leb S\in(0,\infty)$. Let $X$ be a point chosen
uniformly at random from the set $[\Pi]\cap S$ (provided it is
non-empty), conditional on $\Pi$. Then
$$\py\big(X \text{\rm\ exists and is $r$-isolated}\big)
\geq e^{-\lambda\leb B(0,r)}-e^{-\lambda\leb S}.$$
\end{lemma}

\begin{lemma}
\label{correlation}
Under the assumptions of Lemma \ref{isolated}, let $Z\subseteq\rd$
be disjoint from $S$. Then
$$\py\big(X \text{\rm\ exists and is $r$-isolated}\bigmid \Pi(Z)=0\big)
\geq \py\big(X \text{\rm\ exists and is $r$-isolated}\big).$$
\end{lemma}

\begin{pfof}{Lemma \ref{isolated}}
For convenience, let $X=\Delta$ if $\Pi(S)=0$, and take $\Delta$ to be
not $r$-isolated. Denote the random sets $B=B(X,r)$, $U=B\cap S$, and
$V=B\setminus S$. Then we have
\begin{eqnarray*}
\lefteqn{\py(X \text{\rm\ is $r$-isolated}\mid \Pi(S), X)}\\ & =&
\py\big(\Pi(U)=1, \Pi(V)=0 \mid \Pi(S), X\big)\ind[\Pi(S)>0] \\ &=&
\bigg( 1-\frac{\leb U}{\leb S}\bigg)^{\Pi(S)-1} e^{-\lambda\leb V}
\ind[\Pi(S)>0].
\end{eqnarray*}
Hence
\begin{eqnarray*}
\py(X \text{\rm\ is $r$-isolated}) &=& \ey\py(X \text{\rm\ is
$r$-isolated}\mid \Pi(S), X) \\
&=& \ey\sum_{k=1}^\infty
\bigg(1-\frac{\leb U}{\leb S}\bigg)^{k-1} e^{-\lambda\leb V} e^{-\lambda\leb S}
\frac{(\lambda\leb S)^k}{k!} \\
&=& \ey\Bigg[\bigg( 1-\frac{\leb U}{\leb S}\bigg)^{-1}
 e^{-\lambda\leb V-\lambda\leb S}
\big(e^{\lambda\leb S-\lambda\leb U}-1\big)\Bigg] \\
&=& \ey\Bigg[\bigg( 1-\frac{\leb U}{\leb S}\bigg)^{-1} \big(e^{-\lambda\leb
B}-e^{-\lambda\leb V-\lambda\leb S}\big)\Bigg] \\
&\geq & e^{-\lambda\leb B}-e^{-\lambda\leb S}.
\end{eqnarray*}
\end{pfof}

Let $\Pi|_S$ denote $\Pi$ restricted to $S$; that is the point
process with support $[\Pi]\cap S$.

\begin{pfof}{Lemma \ref{correlation}}
We have
\begin{eqnarray*}
\lefteqn{\py\big(X \text{ is $r$-isolated},\Pi(Z)=0\bigmid \Pi|_S,X\big)}\\
&=&\ind\big[\Pi(B(X,R)\cap S)=1\big]
 \py\big(\Pi(B(X,R)\setminus S)=0,\Pi(Z)=0\bigmid X\big) \\
&\geq&\ind\big[\Pi(B(X,R)\cap S)=1\big]\py\big(\Pi(B(X,R)\setminus
S)=0\bigmid X\big)
 \py\big(\Pi(Z)=0\big) \\
&=&\py\big(X \text{ is $r$-isolated}\bigmid
\Pi|_S,X\big)\py\big(\Pi(Z)=0\big).
\end{eqnarray*}
Taking expectations yields the result.
\end{pfof}

\subsection*{Volume bound}

We write $\omega_d:=\leb B_d(0,1)$ for the volume of the unit
ball. Fix $d\geq 3$ and $\epsilon,C>0$. Suppose
$w,w'\in\mathbb{R}^2$ are such that $\rho(w,w')=1$. Define the
sets
$$
W=B_2(w,\epsilon)\times B_{d-2}(0,C) \quad\text{and}\quad
W'=B_2(w',\epsilon)\times B_{d-2}(0,C).
$$

\pagebreak
\begin{lemma}
\label{volcalc} Fix $\mu=0.75$, $\delta=0.1$ and $\epsilon=0.01$.
Let $d\ge 11$.  Let $W,W'$ be as above.  There exists $C'=C'(d)$
such that if $C\geq C'$, then for all $x\in W$ and $r\in [\mu-\delta,
\mu+\delta]$, writing
$$S=W'\cap \big(B(x, r+\mu+\delta)\setminus B(x,r+\mu-\delta)\big),$$
\nopagebreak we have \nopagebreak
$$\leb S\ge\frac{\omega_2\omega_{d-2}{\epsilon}^2}{3}
\left(1.2^{\frac{d-2}{2}}-1\right).$$
\end{lemma}

In order to prove Lemma~\ref{volcalc}, we
use some further geometric results.
\begin{lemma}
\label{twosph} Fix $d$ and $0<R_{-}\leq R_{+}<\infty$. There
exists $C'=C'(d,R_{-},R_{+})$  such that if $C\geq C'$, for all $x\in
B(0,C)$ and $R\in (R_{-},R_{+})$ we have
$$\leb \big(B(0,C)\cap B(x,R)\big)\ge\frac{1}{3}\leb \big(B(x,R)\big).$$
\end{lemma}
\begin{pf}
First suppose that $\rho(0,x)=C$ (so that $x$ is on the surface of
the ball $B(0,C)$). Let
$$f(C,R)=\frac{\leb (B(0,C)\cap B(x,R))}{\leb (B(x,R))},$$
and note that $f$ depends only on the ratio $R/C$, while for fixed
$R$, the function $f$ is increasing and continuous in $C$, and
converges to $1/2$ as $C\to\infty$ (since near $x$, the ball
$B(0,C)$ approaches a half space).  Therefore by the intermediate
value theorem the claimed result holds for the case $\rho(0,x)=C$.

Suppose now that $\rho(0,x)<C$.  The result is trivial when
$B(x,R)\subseteq B(0,C)$.  If not we can replace $B(0,C)$ with
$B(0,\widetilde{C})$, where $\rho(0,x)=\widetilde{C}$ (and so
$\widetilde{C}\in[C-R_+,C)$), and appeal to the case already
proved.
\end{pf}

\begin{lemma}
\label{sphcyl} Let $\mu,\delta,\epsilon$ and $W,W'$ be as in Lemma
\ref{volcalc}. There exists $C'(d)<\infty$ such that if $C\geq C'$,
then for all $x\in W$ and $R\in [2(\mu-\delta), 2(\mu+\delta)]$ we
have
$$\frac{1}{3}{\omega}_2 {\omega}_{d-2} {\epsilon}^2 R_1^{d-2}
\le \leb \big(W'\cap B(x,R)\big) \le {\omega}_2 {\omega}_{d-2}
{\epsilon}^2 R_2^{d-2},$$
where $R_{1}=\sqrt{R^2-(1+2\epsilon)^2}$ and
$R_{2}=\sqrt{R^2-(1-2\epsilon)^2}.$
\end{lemma}
\begin{pf}
For $z=(z_1, z_2,\ldots z_d)\in\rd$ we write $\underline{z}=(z_1,
z_2)$ and $\overline{z}=(z_3, z_4,\ldots z_d)$, and we write
$\underline{\rho}$ and $\overline{\rho}$ for Euclidean distances
on ${\mathbb R}^2$ and ${\mathbb R}^{d-2}$ respectively.

We claim that for any $x\in W$,
\begin{multline*}
B_2(w',\epsilon)\times \big[B_{d-2}(0,C)\cap B_{d-2}(\overline{x},R_1)\big] \\
\subseteq W'\cap B(x,R) \subseteq \\
B_2(w',\epsilon)\times \big[B_{d-2}(0,C)\cap B_{d-2}(\overline{x},R_2)\big].
\end{multline*}
To prove this, first note that for $\underline{x}\in
B_2(w,\epsilon)$ and $\underline{z}\in B_2(w',\epsilon)$, we have
$$\underline{\rho}(w,w')-\underline{\rho}(\underline{z},w')
-\underline{\rho}(\underline{x},w)\le
\underline{\rho}(\underline{z},\underline{x})\le
\underline{\rho}(\underline{z},w')+\underline{\rho}(w,w')+
\underline{\rho}(w',\underline{x}),$$
which gives
$$1-2\epsilon\le \underline{\rho}(\underline{z},\underline{x})\le
1+2\epsilon.$$

Take any $z\in W'\cap B(x,R)$. Then $\underline{z}\in
B_2(w',\epsilon)$ and $\underline{x}\in B_2(w,\epsilon)$.  Since
$z\in B(x,R)$, we have
$$\overline{\rho}(\overline{z},\overline{x})
\le \sqrt{R^2-\underline{\rho}(\underline{z},\underline{x})^2}
\le\sqrt{R^2-(1-2\epsilon)^2}=R_2;$$
that is $\overline{z}\in B_{d-2}(\overline{x},R_2)$.
The second inclusion of the claim follows.

On the other hand, if $z\in B_2(w',\epsilon)\times (B_{d-2}(0,C)\cap
B_{d-2}(\overline{x},R_1))$, then
$\underline{\rho}(\underline{z},\underline{x})\le 1+2\epsilon$ and
$\overline{\rho}(\overline{z},\overline{x})\le R_1.$ This gives
$\rho(z,x)\le R$, so $z\in W'\cap B(x,R).$ Thus the first inclusion of
the claim is proved.

The required statement now follows from the claim.  For the lower
bound we use Lemma \ref{twosph}, noting that
$$R_1\in\big[\sqrt{4(\mu-\delta)^2-(1+2\epsilon)^2}\;,\;
\sqrt{4(\mu+\delta)^2-(1+2\epsilon)^2}\big]\subset(0,\infty).$$
For the upper bound we discard the intersection with $B_{d-2}(0,C)$.
\end{pf}

\begin{pfof}{Lemma \ref{volcalc}}
We have
$$\leb S=\leb \big(W'\cap B(x,r+\mu+\delta)\big)
-\leb \big(W'\cap B(x, r+\mu-\delta)\big).$$
Therefore, choosing $C'$ according to Lemma~\ref{sphcyl},
for $C\geq C'$ we have by that lemma,
$$\leb S \ge \omega_2\omega_{d-2}{\epsilon}^2 \Big[
\overline{r}^{d-2}/3-\underline{r}^{d-2}\Big],$$
where
$\overline{r}:=\sqrt{(r+\mu+\delta)^2-(1+2\epsilon)^2}$ and
$\underline{r}:=\sqrt{(r+\mu-\delta)^2-(1-2\epsilon)^2}.$ For $C$
chosen as above, we want to find a uniform lower bound on $\leb S$ for
all possible $r\in [\mu-\delta,\mu+\delta]$.  We claim that for $d\geq 11$,
the function
$$g(r):=\overline{r}^{d-2}/3-\underline{r}^{d-2}$$
is increasing on $r\in [\mu-\delta,\mu+\delta]$.
Once this is proved we obtain
$$g(r)\geq g(\mu-\delta) \geq (\sqrt{1.2})^{d-2}/3-(\sqrt{0.73})^{11-2}\geq
\tfrac 13 \left(1.2^{\frac{d-2}{2}}-1\right)$$
on this interval, and the result follows.

To prove the above claim note that
$$g'(r)=\frac{d-2}{3}\Big[(r+\mu+\delta)\overline{r}^{d-4}
-3(r+\mu-\delta)\underline{r}^{d-4}\Big],$$
so it is enough to check that
$$\bigg(\frac{\overline{r}^2}{\underline{r}^2}\bigg)^{\frac{d-4}{2}}>3$$
for the required values of $r$ and $d$.
But is is straightforward to check that the quadratic function
$\overline{r}^2-1.4\underline{r}^2$ is positive for $r\in
[\mu-\delta,\mu+\delta]$, hence $\overline{r}^2/ \underline{r}^2>1.4$;
and we have $1.4^{\frac{11-4}{2}}>3$.
\end{pfof}

\subsection*{Proof of Theorem 2}

\begin{pfof}{Theorem \ref{noninv}}
Let $d\geq 45$, and construct the hard sphere process $\Gamma$ as
in Section \ref{construction}, where the constant $C$ is chosen
according to Lemma \ref{volcalc}, and the intensity $\lambda$ of
the Poisson process will be chosen later.

Consider the first layer of the construction, and recall that at
step $n$, vertex $w_n$ is explored, and a sphere of radius
$r_{w_n}$ is placed with centre $x_{w_n}=y_{w_n}\in W(w_n)$, provided the
vertex is found to be good.  As a notational convenience, if the algorithm
stops during step $N$ we write $w_n=\Delta$ for all $n\geq N$, and
call $\Delta$ a good vertex.  Let ${\cal F}_n$ be the
$\sigma$-algebra generated by all of the random variables
$$
\Big(w_i,y_{w_i}, x_{w_i}, r_{w_i} ,\Pi|_{W(w_i)},
\Pi|_{B(x_{w_i},r_{w_i})} \Big)_{i=0,\ldots, n};
$$
that is ``the information known up to and including step $n$''.
Let ${\cal F}_{-1}$ be the trivial $\sigma$-algebra.

We will compare the set of good vertices with a percolation
cluster. Suppose that for some $q$ we can show that for all $n\geq
0$,
\begin{equation}
\label{success} \py(w_n \text{ is good}\mid {\cal{F}}_{n-1})\ge q
\quad\text{almost surely}
\end{equation}
(that is, each newly explored vertex is good with probability at
least $q$ uniformly in the past).  Then the random set of good
{\em site} vertices stochastically dominates
$C_{q^2}^{\mathbb{H}}(0),$ the open cluster at the origin for site
percolation with parameter $q^2$ on ${\mathbb H}$.  This is
because, as long as it is possible to add a new site vertex to the
cluster of good site vertices at the origin, the algorithm
attempts to do so by first exploring the intervening bond vertex,
and then immediately exploring the new site vertex.  The
probability that both steps succeed is at least $q^2$.  Therefore
we can compare with a cluster-growing algorithm for
$C_{q^2}^{\mathbb{H}}(0)$.

Also recall that the construction gives, in each layer, a cluster of
touching spheres (with radii $\in[\mu-\delta, \mu+\delta]$) including
a sphere with its centre in the cell of each good vertex.  Note
therefore that if $J$ is such a cluster of spheres then
$J^{\{C+\epsilon\}}$ contains all good vertices.  Therefore, since the
layers of the construction are independent, if we can establish
(\ref{success}) with $q^2>\pch$ then the statements of the theorem
will follow by Lemma \ref{perc-dens}.  Thus all that remains is to
prove (\ref{success}) for with $q\geq 0.892> \sqrt{0.794}\geq
\sqrt{\pch}$.

Recall that in step $n$ of the construction we choose a random
Poisson point $y=y_{w_n}$ (if any exists) in a certain
${\cal{F}}_{n-1}$-measurable set $S$, and check to see whether a
certain ball $B(y,r_{w_n})$ is free of other Poisson points.  If
both steps succeed then the vertex is declared good.  The radius
$r_{w_n}$ depends on which point $y$ is chosen, but since it can
be at most $\mu+\delta$ we can bound the required probability
by the probability that the
larger ball $B(y,\mu+\delta)$ contains no other points:
$$ \py(w_n \text{ is good}\mid {\cal{F}}_{n-1})\ge
\py\big(y\neq\Delta, \text{ and } y \text{ is
$(\mu+\delta)$-isolated} \bigmid {\cal{F}}_{n-1} \big).$$ The
latter event depends on $\Pi$ only through its restriction to
$W(w_n)^{\{\mu+\delta\}}$.  The conditioning on ${\cal F}_{n-1}$
does not affect the process $\Pi|_{W(w_n)}$, while the fact that
one vertex $v_n$ adjacent to $w_n$ is good tells us only that a
certain ${\cal{F}}_{n-1}$-measurable set $Z$ not intersecting $S$
contains no points of $\Pi$.  (Recall that
$2(\mu+\delta+\epsilon)<\surd 3$,
so the conditioning on non-adjacent vertices has no effect).
Thus we deduce by Lemmas
\ref{correlation} and \ref{isolated} that
\begin{align*}
\py(w_n \text{ is good}\mid {\cal{F}}_{n-1})&\ge
e^{-\lambda \leb B(0,\mu+\delta)}-e^{-\lambda \leb S}\\
&\geq 1-\lambda \leb B(0,\mu+\delta)-e^{-\lambda \leb S}.
\end{align*}
But Lemma \ref{volcalc} applies to give
$$\leb S \geq
\frac{\omega_2\omega_{d-2}{\epsilon}^2}{3}\left(1.2^{\frac{d-2}{2}}-1\right)$$
almost surely.

Thus, we require that
$$F(\lambda):=1-\lambda B-
e^{-\lambda A} \geq 0.892,$$ where
$$
A:=\frac{\omega_2\omega_{d-2}10^{-4}}{3}\left(1.2^{\frac{d-2}{2}}-1\right)
\quad\text{and}\quad
B:=\omega_d 0.85^d.
$$
For each $d$ we can choose the intensity $\lambda$ so as to get
the best bound.  Differentiating shows that $F$ has a maximum at
$$\lambda^*=\lambda^*(d)=\frac1A \log \frac AB,$$
at which
$$F(\lambda^*)= 1-\frac{\log{(A/B)}+1}{A/B}.$$
But $\omega_2=\pi$ and
$\frac{\omega_{d-2}}{\omega_d}=\frac{d}{2\pi}$ (see for example
\cite{somm}), therefore
$$\frac AB =\frac{10^{-4} d \left(1.2^{\frac{d-2}{2}}-1\right)}
{6\times 0.85^d}.$$
Thus $A/B$ is an increasing function of $d$, and it is easy to
check that for $d\geq 31$ we have $A/B>1$, and therefore
$\lambda^*$ is positive (which is a requirement for an
intensity). Furthermore, $F(\lambda^*)$ is an increasing function
of $A/B$ for $A/B>1$, and therefore an increasing function of
$d\geq 31$. Finally it is straightforward to check that
$F(\lambda^*(d))> 0.892$ for $d\geq 45$, as required.
\end{pfof}

\paragraph{Remark -- choice of lattice.}
Our construction could be adapted to work for other
two-dimensional lattices such as the square lattice in place of
the hexagonal lattice, but at the expense of increasing the
dimension $d$.  The fundamental requirement is that the diagonals
(that is, the distances between non-adjacent vertices of the
lattice) must be greater than $2/\surd 3$ times the edge length.
This allows the set $S$ (in which we try to find a possible centre
for a sphere) to be made much larger in volume than the hard
spheres (which must be empty of other Poisson points), while
preventing overlap between non-adjacent spheres.

\section{Stationarization}
\label{average}

In this section we deduce Theorem \ref{main1} from Theorem \ref{noninv}.

\begin{pfof}{Theorem \ref{main1}}
For each positive integer $n$, let $U_n$ be a random variable
uniformly distributed on $B(0,n)$, and independent of $\Gamma$.
Define the randomly shifted process $\Gamma_n=\Gamma+U_n$. Clearly
$\Gamma_n$ is a percolating Poisson hard sphere process, and has
no spheres larger than $K$. We shall use Prohorov's Theorem to
construct $\Lambda$ as a weak limit of the sequence $(\Gamma_n)$,
and show that it has all the required properties.

We claim that the sequence of random variables $(\Gamma_n)$ is tight
in the weak topology induced by the vague topology on point measures
on $\rd\times[0,\infty)$. To check this, it is enough to show that the
sequence $(\Gamma_n(A))$ is tight for any relatively compact Borel
$A\in \rd\times[0,\infty)$. (See \cite{kallen} Lemma 16.15).  Any such
$A$ is a subset of $B\times[0,\infty)$ for some bounded Borel
$B\subseteq\rd$, and we have \as
$\Gamma_n(A)\leq\Gamma_n(B\times[0,K])$ since $\Gamma_n$ has no
spheres larger than $K$.  But the latter quantity has the same
(Poisson) distribution for each $n$, so the sequence is clearly tight
as required.

Now let $\Lambda$ be any weak subsequential limit of $(\Gamma_n)$, so
\begin{equation}
\Gamma_{n_k}\Rightarrow \Lambda \text{ as }k\to\infty
\label{weakconv}
\end{equation}
in the topology referred to above. Clearly $\Lambda$ is
integer-valued and thus a sphere process. Furthermore, it is
easily seen that the set of hard sphere processes supported \as on
$\rd\times[0,K]$ is weak closed, and therefore $\Lambda$ is a hard
sphere process supported on $\rd\times[0,K]$. Let $B\subseteq\rd$
be a bounded Borel set with $\leb$-null boundary. By Theorem 16.16
of \cite{kallen}, the convergence in (\ref{weakconv}) implies the
convergence in distribution $\Gamma_{n_k}(B\times[0,K])\to
\Lambda(B\times[0,K])$. Hence the latter has Poisson distribution
with mean $\leb B$ and so $\Lambda$ is a Poisson sphere process.

Next we show that $\Lambda$ is invariant. It is sufficient to show
that for any $z\in\rd$ and any continuous compactly-supported
function $f:\rd\times[0,\infty)\to [0,1]$ we have $\ey\int f\;
d\Lambda = \ey\int f\; d(\Lambda +z)$. Recall that
$\Gamma_n=\Gamma+U_n$; we shall compare $\Gamma_n$ and
$\Gamma_n+z$ for large $n$. Fix $z$ and $f$, and write
$J=J(n)=B(0,n)\cap B(z,n)$. Let $V$ be random variable uniformly
distributed on $J$. Also write $\alpha=\alpha(n)=\leb
(B(0,n)\setminus J)/\leb B(0,n)$, and note that $\alpha\leq
cn^{d-1}/n^d\to 0$ as $n\to\infty$. Recall that $U_n$ was uniform
on $B(0,n)$. Observe that conditional on $U_n\in J$ (which occurs
with probability $1-\alpha$), the law of $U_n$ equals that of $V$.
And conditional on $U_n+z\in J$ (which occurs with probability
$1-\alpha$), the law of $U_n+z$ equals that of $V$. Hence
\begin{eqnarray*}
\lefteqn{\ey \int f \; d\Gamma_n}\\
 &=& \ey\bigg[\int f \;d(\Gamma+U_n)\bigg| U_n\in
J\bigg](1-\alpha) + \ey\bigg[\int f \;d(\Gamma+U_n)\bigg| U_n\in
J^c\bigg]\alpha \\
&=& \ey\int f \;d(\Gamma+V)(1-\alpha) + \ey\bigg[\int f
\;d(\Gamma+U_n)\bigg| U_n\in J^c\bigg]\alpha,
\end{eqnarray*}
and similarly
\begin{eqnarray*}
\lefteqn{\ey \int f \; d(\Gamma_n+z)} \\
 &=&
\ey\int f \;d(\Gamma+V)(1-\alpha) +
\ey\bigg[\int f \;d(\Gamma+U_n+z)\bigg| U_n+z\in
J^c\bigg]\alpha.
\end{eqnarray*}
Since $U_n$ is independent of $\Gamma$, the processes $\Gamma+U_n$ and
$\Gamma+U_n+z$ appearing on the right sides of the above equations are
Poisson sphere processes, even when conditioned on $U_n$.
Now for any Poisson sphere process $\Upsilon$ say, $\ey\int f d\Upsilon$
is bounded by the expected number of Poisson points in the projection
of the support of
$f$ onto $\rd$, that is the Lebesgue measure of that projection, $C$ say. Hence,
subtracting the two equations above gives
$$\left|\ey\int f\;d\Gamma_n - \ey\int f\;d(\Gamma_n+z)\right|\leq
0+ \alpha C \to 0
\text{ as } n\to\infty.$$
Taking weak limits as $k\to\infty$ of $\Gamma_{n_k},\Gamma_{n_k}+z$ we
deduce from (\ref{weakconv}) that
$$\left|\ey\int f\;d\Lambda - \ey\int f\;d(\Lambda+z)\right|=0.$$
Hence $\Lambda$ is invariant as required.

Finally we must show that $\Lambda$ percolates almost surely.
 For a hard sphere process $\Upsilon$ and for
$0<a<b$ let $H_{a,b}(\Upsilon)$ be the event that $G(\Upsilon)$
has a connected set of spheres with radii at most $K$ which
intersects both $B(0,a)$ and $B(0,b)^C$. Also let
$H_{a,\infty}(\Upsilon)$ be the event that $G(\Upsilon)$ has an
unbounded connected set of spheres with radii at most $K$ which
intersects $B(0,a)$. Note that $H_{a,\infty}(\Upsilon)$ is the
decreasing limit of the events $H_{a,b}(\Upsilon)$ as
$b\to\infty$. Also denote by $I(\Upsilon)$ the union of all
infinite clusters of the hard sphere process $\Upsilon$. Recalling
the definition of $\Gamma_n$ above, we have
\begin{eqnarray*}
\py(H_{a,b}(\Gamma_n))&\geq&\py(H_{a,\infty}(\Gamma_n))\\
&=&\py(0\in I(\Gamma)^{\{a\}})\\
&=&\ey\frac{\leb(I(\Gamma)^{\{a\}}\cap B(0,n))}{\leb B(0,n)}.
\end{eqnarray*}
Hence by Fatou's Lemma and the definition of lower density,
$$\liminf_{n\to\infty}\py(H_{a,b}(\Gamma_n))\geq\ey\ld(I(\Gamma)^{\{a\}}).$$
Note that the event $H_{a,b}(\Upsilon)$ depends only on the
process $\Upsilon$ restricted to the compact set
$B(0,b+K)\times[0,K]$. Furthermore it is straightforward to see
that the event is closed in the vague topology on point measures.
It follows from (\ref{weakconv}) and the Portmanteau Theorem
(\cite{kallen} Theorem 4.25) that
$$\py(H_{a,b}(\Lambda))\geq\limsup_{k\to\infty}\py(H_{a,b}(\Gamma_{n_k})).$$
From the last two inequalities we have
$$\py(H_{a,b}(\Lambda))\geq \ey\ld(I(\Gamma)^{\{a\}}),$$
and hence letting $b\to\infty$ we deduce
$$\py(\Lambda \text{ percolates})\geq \py(H_{a,\infty}(\Lambda))\geq
 \ey\ld(I(\Gamma)^{\{a\}}).$$
Finally letting $a\to\infty$, Theorem \ref{noninv} gives
$$\py(\Lambda \text{ percolates})=1.$$
\end{pfof}

\paragraph{Remark -- positive radii.}
As noted earlier, our proof may be adapted so that all the spheres
have positive radii.  To achieve this, take a small parameter
$\eta>0$, and modify the construction in Section
\ref{construction} as follows: having chosen a potential centre
$y_w$ for a sphere, we declare the vertex good only if the larger
ball $B(y_w,r_w+\eta)$ contains no other Poisson points (rather
than the ball $B(y_w,r_w)$). If $\eta$ is small enough then this
does not affect the later computations, and we still obtain a
non-invariant percolating hard sphere process $\Gamma$ for $d\geq
45$.  But now $\Gamma$ has the additional property that no
zero-radius sphere is within distance $\eta$ of any nonzero-radius
sphere.  Therefore the stationarized version $\Lambda$ will
inherit the same property.  Finally, we modify $\Lambda$ as
follows.  If there is a zero-radius sphere centred at $z$, replace
it with a sphere of radius $r$, where $r$ is $1/2$ of the distance
from $z$ to the closest other sphere of $\Lambda$ (including other
zero-radius spheres).  This $r$ is always positive because the
radii of the existing spheres are uniformly bounded above, and any
bounded region of $\rd$ contains only finitely many Poisson
points.

\section*{Open Problems}

\begin{mylist}
\item Does there exist a percolating Poisson hard sphere process
(invariant or non-invariant) in dimensions $2\leq d\leq 44$?  The
case $d=2$ seems particularly interesting.

\item In any dimension, does there exist a percolating, invariant
Poisson hard sphere process which is a deterministic function of
the Poisson process?

\item Do percolating hard sphere processes exist for other point
processes, such as Gaussian zeros processes \cite{zeros}.
\end{mylist}

\section*{Acknowledgements}

We thank Yuval Peres for drawing our attention to the problem.  Codina
Cotar thanks her postdoctoral advisor David Brydges for support and
assistance.

\section*{ }

Codina Cotar: {\tt c.cotar@math.ubc.ca} \\
Alexander E. Holroyd: {\tt holroyd@math.ubc.ca} \\
University of British Columbia, \\
121-1984 Mathematics Rd, \\
Vancouver BC V6T 1Z2, Canada.

\vspace{4mm} \noindent David Revelle: {\tt
david.revelle@gmail.com}

\end{document}